\documentclass[11pt,letterpaper]{IEEEtran}

  \onecolumn
\usepackage{amsmath, amssymb,amsthm,cite}
\usepackage[utf8x]{inputenc}
\usepackage{graphicx}
\usepackage{psfrag}
\usepackage{bbm}
\usepackage{mathtools}
\usepackage{color}
\usepackage{tcolorbox}
\usepackage{arydshln}

\DeclareMathOperator*{\nn}{\nonumber}
\DeclareMathOperator{\E}{\mathbb{E}}

\DeclareMathOperator*{\regret}{regret}

\allowdisplaybreaks
\newtheorem{lemma}{Lemma}
\newtheorem{problem}{Problem}
\newtheorem{theorem}{Theorem}

\newtheorem{remark}{Remark}
\theoremstyle{definition}

\renewcommand{\v}[1]{\mathbf{#1}}
\renewcommand{\c}[1]{(zI-#1)^{-1}}
\renewcommand{\a}[1]{(z^{-1}I-#1)^{-1}}

\DeclareMathOperator{\vw}{\mathbf{w}}

\DeclareMathOperator{\vv}{\mathbf{v}}
\DeclareMathOperator{\vy}{\mathbf{y}}

\makeatletter
\def\blfootnote{\gdef\@thefnmark{}\@footnotetext}
\makeatother

\DeclarePairedDelimiterX{\norm}[1]{\lVert}{\rVert}{#1}

\title{Regret-Optimal Filtering for Prediction and Estimation}
\author{Oron Sabag and Babak Hassibi}
\begin{document}
\maketitle

\begin{abstract}
    In this paper, we study the filtering problem of causally estimating a desired signal from a related observation signal, through the lens of regret optimization. Classical filter designs, such as $\mathcal H_2$ (i.e., Kalman) and $\mathcal H_\infty$, minimize the average and worst-case estimation errors, respectively. As a result $\mathcal H_2$ filters are sensitive to inaccuracies in the underlying statistical model, and $\mathcal H_\infty$ filters are overly conservative since they safeguard against the worst-case scenario. We propose instead to minimize the \emph{regret} in order to design filters that perform well in different noise regimes by comparing their performance with that of a clairvoyant filter. More explicitly, we minimize the largest deviation of the squared estimation error of a causal filter from that of a non-causal filter that has access to future observations. In this sense, the regret-optimal filter will have the best competitive performance with respect to the non-causal benchmark filter no matter what the true signal and the observation process are. For the important case of signals that can be described with a time-invariant state-space, we provide an explicit construction for the regret optimal filter in the estimation (causal) and the prediction (strictly-causal) regimes. These solutions are obtained by reducing the regret filtering problem to a Nehari problem, i.e., approximating a non-causal operator by a causal one in spectral norm. The regret-optimal filters bear some resemblance to Kalman and $\mathcal H_\infty$ filters: they are expressed as state-space models, inherit the finite dimension of the original state-space, and their solutions require solving algebraic Riccati equations. \textcolor{black}{Numerical simulations demonstrate that regret minimization inherently interpolates between the performances of the $\mathcal H_2$ and $\mathcal H_\infty$ filters and is thus a viable approach for filter design.}
\end{abstract}

\section{Introduction}
\blfootnote{The authors are with the
Department of Electrical Engineering at California Institute of Technology (e-mails:
 \{oron,hassibi\}@caltech.edu). A preliminary version of the results presented here first appeared in \cite{SabagFilteringAISTATS}.}
Filtering is the problem of estimating the current value of a desired signal, given past values of a related observation signal. It has numerous applications in signal processing, control, and learning and its rich history going back to Wiener, Kolmogorov, and Kalman \cite{Wiener_filter_book,KalmanFiltering}. When the underlying desired and observation signals have state-space structures driven by white Gaussian noise, the celebrated Kalman filter gives the minimum mean-square error estimate of the current value of the desired signal, given the past and current of the observed signal. When all that is known of the noise sources are their first and second-order statistics, the Kalman filter gives the linear optimal least-mean-square-estimate. While these are very desired properties, the Kalman filter is predicated on knowing the underlying statistics and distributions of the signal. It can therefore have poor performance if the underlying signals have statistics and/or properties that deviate from those that are assumed. It is also not suitable for learning applications, since it has no possibility of \emph{learning} the signal statistics.

Another approach to filtering that was developed in the 80's and 90's was $\mathcal H_\infty$ filtering, where the noise sources were considered adversarial and the {\em worst-case} estimation error energy is minimized over all bounded-energy noises \cite{grimble1990h,li1997linear,ShakedSurveyEstimation,57538,67291,TamerBasar_book}. While $\mathcal H_\infty$ estimators are robust to lack of statistical knowledge of the underlying noise sources, and have some deep connections to classical learning algorithms (see, e.g. \cite{HSK95}), they are often too conservative since they safeguard against the worst-case and do not exploit the noise structure. We propose to adopt the regret so as to bridge between the philosophies of Kalman ($\mathcal H_2$) and $\mathcal H_\infty$ filtering.

\textcolor{black}{Regret is a competitive criterion that measures the loss in performance of a practical policy with respect to a superior policy that cannot be implemented in practice. In the context of filtering, we formulate a regret criterion as the absolute difference between the squared estimation error of a causal filter, to be designed, and that of a clairvoyant filter that has access to the entire observations sequence (including future observations). We then define the optimal regret as the largest absolute squared estimation error difference across all noise sequences. By minimizing the regret, the designed causal filter aims to mimic the behavior of the non-causal filter across all noise sequences. In other words, regret minimization guarantees that the estimation error of the regret-optimal filter deviates the least from the performance of the non-causal filter that serves as our benchmark. The simulations demonstrate that regret minimization leads to filters that interpolate between the $\mathcal H_2$ and the $\mathcal H_\infty$ approaches to obtain filters that have balanced behavior across all noise sequences.}

The regret has been utilized in wide range of problems in learning, estimation \cite{eldar_regret,eldar_merhav_regret1,eldar_merhav_regret2,yasini2017worst,feder2021sequential} and control \cite{cohen2019learning,muthirayan2021online,hazan2019nonstochastic,simchowitz2020improper,agarwal2019online,cassel2021online,chang2021regret,SabagFIArxiv}. In particular, the regret approach proposed in \cite{SabagFIArxiv} for control problems is closely related to the one taken here. In learning, for instance, the regret is defined as the difference between the cost of a feasible policy that does not have access to the system parameters, and it is compared with the cost achieved by a supreme policy which has a full knowledge of the system parameters. The objective there is to derive a policy that will obtain a \emph{sub-linear regret} meaning that, asymptotically, the inferior policy will attain the optimal (averaged) cost attained by the supreme policy. In this work, we adopt the appealing comparative nature of the regret, but utilize it in a different way than the typical ones used in learning theory; our regret cannot obey a sub-linear law since the causal estimator is inferior with respect to the optimal non-causal estimator no matter how large is the horizon. This fact allows a meaningful comparison with a benchmark . Indeed, the comparison with the non-causal filter enables us to solve the problem in terms of the optimal regret computation and the regret-optimal filter.

Our main results are for the important case of linear filtering problems that can be described with a state-space model. For time-invariant state-space models, we derive regret-optimal filters for two scenarios. \textcolor{black}{The first scenario is an estimation regime where the designed filter constructs its current estimate based on the past and present observations (i.e., it has causal access to the observations), and the other scenario is the prediction regime where a predicted estimate of the current state is constructed based on past observations only (i.e. it has a strictly-causal access to the observations).} The implementation of the regret-optimal filters is given via two simple steps; first, a bisection method with a simple condition is applied to find the optimal regret value. This is done once, prior to the filtering process itself. Then, by having the optimal regret value, we provide a regret-optimal filter in a state-space form. For a state-space with a state dimension $n$, the regret-optimal estimator has $3n$ states in its state-space representation.

The regret problem of comparing \textcolor{black}{the performance} of causal and non-causal policies is closely related with the classical Nehari problem in operator theory of approximating an anti-causal operator with a causal one in an operator norm sense \cite{nehari1957bounded}. We show for general linear operators that the regret-optimal filtering problem can be reduced to the Nehari problem in Theorem \ref{th:reduction}. The reduction to the Nehari problem is done at the operators level and the equivalence is shown by performing several factorizations of positive operators and operators decomposition into their causal and anticausal parts. While the reduction at the operator level does not provide explicit means to construct filters, it holds in the general scenario and sheds light on the structure of regret-optimal filters. For the state-space model, the reduction to the Nehari problem can be made explicit via spectral factorizations and decompositions in order to obtain a Nehari problem in a state-space form, which in turn leads to the explicit regret-optimal predictor and estimator.

We present numerical examples that demonstrate the efficacy and applicability of the regret minimization approach. By comparing the frequency response of the regret-optimal estimator, the $\mathcal H_2$ and the $\mathcal H_\infty$ optimal filters, we observe that the regret-optimal estimator interpolates between the $\mathcal H_2$ and the $\mathcal H_\infty$ performance. We also perform time-domain simulations for different mixtures of white and adversarial noise sequences. In all simulations, it is being observed that the regret-optimal filter has average and worst-case performances that are simultaneously close to their optimal values. Therefore we argue that regret-optimality is a viable approach to filters design. Lastly, it is interesting to note that the regret-optimal estimator maintains an almost constant regret across all frequencies although its performance is determined by the peak of the regret.

The remainder of the paper is organized as follows. Section \ref{sec:formulation} includes the setting and the formulation of the regret approach for the filtering problem. In Section \ref{sec:main}, we present our main results including the reduction to the Nehari problem and the explicit regret-optimal predictor and estimator. Section \ref{sec:examples} includes our numerical experiments. The proof of the reduction of the regret problem to a Nehari problem appears in Section \ref{sec:proof_reduction}, while the state-space derivations are presented in Section \ref{sec:proof}. Technical proofs appear in the appendix.

\section{The Setting and Problem Formulation}\label{sec:formulation}
\textcolor{black}{In this section we present the notation, the problem formulation in a general operator theory notation, the special case of a state-space model, the Riccati and Lyapunov equations needed for the state-space solutions and, lastly, the Nehari problem needed to solve the regret problem.}
\subsection{Notation}
Linear operators are denoted by calligraphic letters, e.g., $\mathcal X$. We use $\{\cdot\}_{+}$ and $\{\cdot\}_{-}$ to denote the causal (lower triangular) and the strictly anticausal (strictly upper triangular) parts of linear operators, respectively. Finite-dimensional vectors and matrices are denoted with small and capital letters, respectively, e.g., $x$ and $X$, and $x^\ast$ denotes the conjugate transpose of $x$. If a matrix $X$ has only real eigenvalues, the maximal eigenvalue is denoted by $\lambda_{\max}(X)$. Subscripts denote time indices e.g., $x_i$, and boldface letters denote the vector consists of finite-dimensional vectors at all times, e.g., $\mathbf{x} = \{x_i\}_i$.
\subsection{The setting and problem formulation}
A linear estimation problem is given by
\begin{align}\label{eq:general_LDS}
    \vy &= \mathcal H \vw + \vv\nn\\
    \v{s} &= \mathcal L \vw
\end{align}
where the sequence $\vw$ denotes an exogenous disturbance, and $\vv$ corresponds to the observations' noise. The sequence $\vy$ denotes the observations process which is used to estimate the desired signal $\v{s}$. The operators $\mathcal H$ and $\mathcal L$ are strictly causal (lower-triangular), and specify the mapping from the disturbance $\vw$ to the observation and the signal processes, respectively. The estimation problem in \eqref{eq:general_LDS} holds for any time horizon and arbitrary linear operators and include, for instance, \textcolor{black}{the state-space model
\begin{align}\label{eq:def_SS_1}
    x_{i+1} &= F x_i + G w_i\nn\\
    y_i&= Hx_i + v_i\nn\\
    s_i &= Lx_i,
\end{align}
that is elaborated in the next section.}

A linear estimator is a mapping from the observations to the space of the signal $\v{s}$ denoted as $\hat{\v{s}} = \mathcal K \vy$. A practical estimator $\mathcal K$ should have a lower-triangular structure to impose the fact that the estimates are causal functions of the observations process. In the scenario of strictly causal filters, which we refer to as the prediction scenario, the operator $\mathcal K$ will have a strictly lower-triangular structure. For both scenarios, the estimation error of a mapping $\mathcal K$ can be written concisely as
\begin{align}\label{eq:TK_def}
    \v{e}(\vw,\vv,\mathcal K)&\triangleq \v{s}-\hat{\v{s}}\nn\\
    &= \mathcal L \vw - \mathcal K [\mathcal H \vw+ \vv]\nn\\
    &= \begin{pmatrix}
    \mathcal L - \mathcal K \mathcal H & -\mathcal K
    \end{pmatrix} \begin{pmatrix}\vw \\ \vv\end{pmatrix}\nn\\
    &\triangleq T_{\mathcal K}  \begin{pmatrix}\vw \\ \vv\end{pmatrix}.
\end{align}
Note that the estimation error is a function of the disturbance $\vw$ and the observation noise $\vv$. The squared estimation error can be also simply expressed as
\begin{align}\label{eq:error_def}
   \v{e}^\ast(\vw,\vv,\mathcal K)\v{e}(\vw,\vv,\mathcal K)&=\begin{pmatrix} \vw^\ast &\vv^\ast \end{pmatrix} T_\mathcal{K}^\ast T_\mathcal{K} \begin{pmatrix}\vw \\ \vv\end{pmatrix}.
\end{align}
The operator $T_\mathcal{K}$ is useful to examine the difference between the existing filtering methods in the literature. For instance,in $\mathcal H_2$ filtering which is known as the stochastic (Kalman) setting where the sequences $\vw$ and $\vv$ are i.i.d. with identity-variance, the expected value of the squared error can be written in terms of the Frobenius norm of the operator $T_\mathcal{K}$, that is, $\E[ \v{e}(\vw,\vv,\mathcal{K})] = \|T_\mathcal{K}\|_F^2$. In the deterministic setting of $\mathcal H_\infty$ filtering, the optimization problem aims to minimize the operator norm of $T_{\mathcal K}$ denoted here as $\|T_{\mathcal K}\|$. A common characteristic of these two paradigms is that if the estimator $\mathcal K$ is not restricted to be causal, there exists a linear estimator that simultaneously attains the minimal Frobenius and operator norms. This fact is summarized in the following lemma.
\begin{lemma}[The non-causal estimator]\label{lemma:noncausal}
For the $\mathcal H_2$ and the $\mathcal H_\infty$ problems, the optimal non-causal estimator is
\begin{align}
 \mathcal{K}_0&= \mathcal{L}\mathcal{H}^\ast (I + \mathcal{H}\mathcal{H}^\ast)^{-1}.
\end{align}
\end{lemma}
The optimality of the non-causal estimator follows from a standard completion of the square and appears, for instance, in \cite[Chapter 10]{hassibi1999indefinite}. Note that the non-causal estimator cannot be implemented in practice even for simple operator $\mathcal L, \mathcal H$ since it requires access to future instances of the observations. The proof of Lemma \ref{lemma:noncausal} can be generalized for norms that are induced from inner products. Lemma \ref{lemma:noncausal} provides a solid motivation to our regret definition since the non-causal estimator outperforms (in the norm sense) any causal or non-causal estimators.

The regret measures the performance of the causal (or strictly causal) filter by computing its maximal error deviation from that of the non-causal filter in Lemma \ref{lemma:noncausal}. In particular, we define the optimal regret as
\begin{align}\label{eq:regret_def}
    \regret&= \min_{\text{causal}\ \mathcal K } \max_{\vw,\vv\neq 0} \frac{\left|\|\v{e}(\vw,\vv,\mathcal K) \|^2 - \|\v{e}(\vw,\vv,\mathcal K_0) \|^2\right|}{\|\vw\|^2 + \|\vv\|^2} \nn\\
    &= \min_{\text{causal}\ \mathcal K } \max_{\vw,\vv\neq 0} \frac{\left|\begin{pmatrix} \vw^\ast &\vv^\ast \end{pmatrix} (T_\mathcal{K}^\ast T_\mathcal{K} - T_{\mathcal{K}_0}^\ast T_{\mathcal{K}_0})\begin{pmatrix}\vw \\ \vv\end{pmatrix}\right|}{\|\vw\|^2 + \|\vv\|^2} \nn\\
    &=  \min_{\text{causal}\ \mathcal K } \|T_\mathcal{K}^\ast T_\mathcal{K} - T_{\mathcal{K}_0}^\ast T_{ \mathcal{K}_0}\|.
\end{align}
\textcolor{black}{In words, the regret measures the maximal deviation of the squared estimation error of the causal from the squared estimation error of the optimal non-causal estimator across all energy-bounded noise sequences. The last equality in \eqref{eq:regret_def} follows from the  maximal absolute Rayleigh quotient of the Hermitian operator $T_\mathcal{K}^\ast T_\mathcal{K} - T_{\mathcal{K}_0}^\ast T_{ \mathcal{K}_0}$.}

It is illuminating to compare the regret with the traditional $\mathcal H_\infty$ estimation:
\begin{align}
    \underbrace{\inf_{\mbox{causal $\mathcal K$}} \|T_{\mathcal K}^*T_{\mathcal K}\|}_{\mbox{$\mathcal H_\infty$ estimation}} ~~,~~
    \underbrace{\inf_{\mbox{causal $\mathcal K$}} \|T_{\mathcal K}^\ast T_{\mathcal K} - T_{\mathcal K_0}^\ast T_{\mathcal K_0}\|.}_{\mbox{regret-optimal estimation}}
\end{align}
The difference is transparent; in $\mathcal H_\infty$ estimation, one attempts to minimize the worst-case gain from the disturbances energy to the estimation error, whereas in regret-optimal estimation one attempts to minimize the worst-case gain from the disturbance energy to the regret. \textcolor{black}{This latter fact implies that the regret-optimal estimator attempts to track, as best as it can, the performance of the non-causal estimator across all possible disturbances. The $\mathcal H_\infty$ estimator, on the other hand, has no no baseline against which to compare itself.} This competitive feature is illustrated in Section \ref{sec:examples}, where our simulations show that regret minimization results in filters that effectively interpolate between the $\mathcal H_2$ and the $\mathcal H_\infty$ design criteria.

Solving the regret-optimal problem in terms of the optimal regret value and the optimal filter is difficult. Thus, we define a sub-optimal regret problem that can be utilized in order to obtain the optimal solution in practice. This is made precise as follows.
\begin{problem}[The $\gamma$-optimal regret estimation problem]\label{pro:regret}
For a fixed $\gamma$, if exists, find a causal (or strictly-causal) $\mathcal K$ such that
\begin{align}
    &\|T_\mathcal{K}^\ast T_\mathcal{K} - T_{\mathcal{K}_0}^\ast T_{\mathcal{K}_0}\| \le \gamma^2.
\end{align}
\end{problem}
A \emph{$\gamma$-optimal estimator (or predictor)} is referred to as any solution to Problem \ref{pro:regret}. The regret-optimal value can be computed by minimizing the value of $\gamma$ for which a $\gamma$-optimal estimator/predictor exists.

\subsection{The state-space model}
The setting in \eqref{eq:general_LDS} is general and it is difficult to obtain an explicit, implementable solution. In many cases, such as Kalman filtering, $\mathcal H_\infty$ estimation, and also the regret problem in Problem \ref{pro:regret}, imposing a state space structure on the operators $\mathcal L$ and $\mathcal H$ can lead to explicit solutions. As mentioned above, in the state-space setting, the equations in \eqref{eq:general_LDS} can be written as
\begin{align}\label{eq:def_SS}
    x_{i+1} &= F x_i + G w_i\nn\\
    y_i&= Hx_i + v_i\nn\\
    s_i &= Lx_i,
\end{align}
where $x_i$ is a hidden state, $y_i$ is the observation and $s_i$ is the signal to be estimated. To recover the state-space setting from its operator notation counterpart in \eqref{eq:general_LDS}, we can choose $\mathcal H$ and $\mathcal L$ as Toeplitz operators with Markov parameters $\mathcal H_i = HF^{i-1}G$ and $\mathcal L_i = LF^{i-1}G$ for $i>0$ and zero otherwise. It is also assumed that $(F,H)$ is detectable and $(F,G)$ is controllable on the unit circle.

A causal estimator is a sequence of mappings $\pi_i(\cdot)$ where the $i$th estimation is $\hat{s}_i=\pi_i(\{y_j\}_{j\le i})$. This means that a causal estimator has access to the past and current observations when constructing its estimate of the current target signal. The estimation error at time $i$ is
\begin{align}
    e_i&= s_i - \hat{s}_i.
\end{align}
In many practical situations, estimators may only have access to past observations, but not to the current observation, when constructing its current estimate. This is a prediction problem and we define a strictly causal estimator, i.e. a predictor, as a sequence of strictly causal mappings, i.e., $\hat{s}_i = \pi_i(\{y_j\}_{j<i})$. Also in the state-space setting, we restrict our attention to linear policies.

Note that the time horizon of the problem may be finite, one-sided infinite and doubly-infinite. However, in the state-space scenario, we focus on the doubly-infinite time horizon regime where the overall estimation error energy is $\sum_{i=-\infty}^\infty  e_i^\ast e_i$. In this horizon regime, it is also convenient to define the transfer matrices
\begin{align}
    H(z)&= H\c{F}G,\nn\\
    L(z)&= L\c{F}G.
\end{align}
to describe the causal mappings from the disturbance $w$ to the observation process $y$ and the target signal $s$, respectively.

\subsection{Riccati and Lyapunov equations}\label{subsec:riccati}
Throughout the derivations, there are three Riccati equations and a single Lyapunov equation. The first Riccati equation is the standard Kalman filter Riccati equation, i.e.,
\begin{align}\label{eq:Riccati_P}
    P\mspace{-2mu} &= \mspace{-2mu}GG^\ast \mspace{-4mu} + \mspace{-2mu} FPF^\ast \mspace{-4mu}- \mspace{-2mu}FPH^\ast (I \mspace{-2mu}+\mspace{-2mu} HPH^\ast)^{-1} HPF^\ast.
\end{align}
The stabilizing solution is denoted as $P$, its feedback gain as $K_P=  FP H^\ast (I + HPH^\ast)^{-1}$ and its closed-loop system as $F_P = F - K_PH$.

The other two Riccati equations depend on the parameter $\gamma$ of Problem \ref{pro:regret}. With some abuse of notation, we omit the dependence on $\gamma$ and will elaborate on it after Theorem \ref{th:strictly_condition}. Define the $\gamma$-dependent Riccati equations as
\begin{align}\label{eq:def_Riccati_gamma}
  W &= H^\ast H + \gamma^{-2}L^\ast L + F^\ast W F - K_W^\ast R_W K_W\nn\\
  Q &= - G R_W^{-1}  G^\ast + F_W Q F_W^\ast  -  K_Q R_Q K_Q^\ast,
\end{align}
with
\begin{align}
 K_W &= R_W^{-1}G^\ast W F\nn\\
 R_W &= I + G^\ast W G \nn\\
 F_W &= F - GK_W
\end{align}
for the first equation, and
\begin{align}
 K_Q &= F_W Q L^\ast R_Q^{-1}\nn\\
 R_Q &= \gamma^2I + L Q L^\ast\nn\\
 F_Q &= F_W-K_QL
\end{align}
for the second equation. Additionally, we define the factorizations $R_W = R_W^{\ast/2}R_W^{1/2}$ and $R_Q = R_Q^{1/2}R_Q^{\ast/2}$. Note that the Riccati equation with $Q$ depends on $W$, the solution to the first Riccati equation. Lastly, define $U$ as the solution to the Lyapunov equation
\begin{align}\label{eq:Lya_U}
 U &= K_QL PF_P^\ast + F_QU F_P^\ast.
\end{align}
We remark that the solutions for \eqref{eq:def_Riccati_gamma} and \eqref{eq:Lya_U} can be obtained with standard software. An implementation of the results in this paper are available at \cite{github_filtering}.

\subsection{\textcolor{black}{The Nehari problem}}\label{subsec:nehari}
Before proceeding to the main results, we define a fundamental problem and its solution that serve as an important tool in our derivations.
\begin{problem}[The Nehari problem]\label{pro:nehari}
Given a strictly anticausal (or anticausal), bounded operator $\mathcal U$, find a causal (or strictly) operator $\mathcal K$ such that
\begin{align}\label{eq:def_nehari_problem}
    &\|\mathcal K - \mathcal U\|
\end{align}
is minimized.
\end{problem}
This problem is known as the Nehari problem \cite{nehari1957bounded}. In the general operator notation, it is difficult to derive an explicit formulae for the approximation $\mathcal K$ and the minimal value of in \eqref{eq:def_nehari_problem}. However, when there is a state-space structure to the strictly anticausal operator $\mathcal U$, the problem has a closed-form solution. The Nehari problem and its solution are known in several cases \cite{Ball1990,GloverDoyle_book}. Here, we use an explicit solution to the Nehari problem for the discrete-time state-space in \cite{SabagFIArxiv}.
\begin{theorem}[Solution to the Nehari problem with a state-space model\cite{SabagFIArxiv}]\label{th:Nehari_general}
Consider the Nehari problem with $T(z) = H (z^{-1}I-F)^{-1}G$ in a minimal form and stable $F$. The optimal norm is given by
\begin{align}\label{eq:th_general_nehari}
    \min_{\mbox{causal, bounded\ $L(z)$}} \| L(z)- T(z) \|&= \lambda_{\text{max}}(Z\Pi),
\end{align}
where $Z$ and $\Pi$ are the unique solutions to the Lyapunov equations
\begin{align}
 Z &= F^\ast Z F + H^\ast H\nn\\
\Pi &= F \Pi F^\ast + G G^\ast.
\end{align}
Moreover, the optimal solution to \eqref{eq:th_general_nehari} is given by
\begin{align}
    L(z)& = H \Pi (I + F_\gamma(zI- F_\gamma)^{-1} )K_\gamma,
\end{align}
with
\begin{align}
    K_\gamma &=  ( I - F^\ast Z_\gamma F \Pi  )^{-1}F^\ast Z_\gamma G\nn\\
    F_\gamma &= F^\ast - K_\gamma G^\ast,
\end{align}
and $Z_\gamma$ is the solution to the Lyapunov equation
\begin{align}
     Z_\gamma &= F^\ast Z_\gamma F + \lambda_{\text{max}}(Z\Pi)^{-1}H^\ast H.
\end{align}
\end{theorem}

\section{Main results}\label{sec:main}
In this section, the main results are presented. In particular, explicit regret-optimal filters for the state-space setting in the strictly-causal and the causal scenarios are presented in Section \ref{subsec:main_strictly} and Section \ref{subsec:main_causal}, respectively. We start with the reduction of the regret problem to the Nehari problem in the estimation problem (causal scenario).
\begin{theorem}[Reduction to the Nehari problem]\label{th:reduction}
A $\gamma$-optimal estimator for the regret problem in the causal scenario (Problem \ref{pro:regret}) exists if and only if the solution to the Nehari problem satisfies
\begin{align}\label{eq:th_reduction_nehari}
\min_{\text{causal}\  \mathcal{K}'}\| \{\nabla_\gamma \mathcal K_0 \Delta \}_- - \mathcal{K}'\|\le 1,
\end{align}
where $\{\cdot\}_-$ denotes the strictly anticausal part of its argument, $\mathcal K_0$ is the non-causal estimator in Lemma \ref{lemma:noncausal}, and $\Delta,\nabla_\gamma$ are causal operators obtained from the canonical factorizations
\begin{align}\label{eq:th_op_factorizations}
    \Delta\Delta^\ast &= I + \mathcal H \mathcal H^\ast\nn\\
    \nabla_\gamma ^\ast\nabla_\gamma &= \gamma^{-2}(I + \gamma^{-2}\mathcal L(I + \mathcal H^\ast \mathcal H)^{-1} \mathcal L^\ast ).
\end{align}

Let $(\gamma^\ast,\mathcal{K}')$ be a solution that achieves the upper bound in \eqref{eq:th_reduction_nehari}, then a regret-optimal causal estimator is given by
\begin{align}\label{eq:th_reduction_filter}
    \mathcal{K} &= \nabla^{-1}_{\gamma^\ast}(\mathcal{K}' + \{\nabla_{\gamma^\ast} \mathcal K_0 \Delta \}_+)\Delta^{-1}
\end{align}
where $\{\cdot\}_+$ denotes the causal part of an operator.
\end{theorem}
Theorem \ref{th:reduction} provides the outline for the solution of the regret problem. The first step is to find the minimal $\gamma$ such that \eqref{eq:th_reduction_nehari} is satisfied. This step concludes the computation of the optimal regret value. Given the optimal $\gamma$ (and the regret value), the next step is to solve the Nehari problem, i.e., to find the causal operator $\mathcal K'$ that solves \eqref{eq:th_reduction_nehari}. This operator leads in turn to the characterization of the regret-optimal filter in \eqref{eq:th_reduction_filter}. The above procedure relies on the canonical factorizations in \eqref{eq:th_op_factorizations}; in the finite-horizon regime, their computation simplifies to Cholesky factorizations but, in the infinite horizon regime of our interest, the computation of factorizations is not trivial. In the next section, we show that the construction of regret-optimal filters is simple when the operators $\mathcal L$ and $\mathcal H$ have a state-space structure. More specifically, frequency domain methods are utilized to compute the factorizations and the decomposition explicitly, and to derive an explicit solution to the Nehari problem in \eqref{eq:th_reduction_nehari}. The proof of Theorem \ref{th:reduction} appears in Section~\ref{sec:proof_reduction}.

\begin{remark}\label{remark}
Theorem \ref{th:reduction} is presented for the causal scenario but can be directly adapted for the strictly-causal scenario by modifying the notation $\{\cdot\}_{-}$ to denote the anticausal part of an operator, and $\{\cdot\}_{+}$ to denote the strictly causal part of an operator. The minimization domain in \eqref{eq:th_reduction_nehari} should be changed accordingly to strictly causal $\mathcal K'$. The resulted Nehari problem is slightly different than \eqref{eq:th_reduction_nehari} since the approximation now is of an anticausal operator with a strictly causal one. Nevertheless, in Section \ref{subsec:proof_str_causal}, we show that these approximation problems can be transformed to one another when the operators have a state-space structure.
\end{remark}

\subsection{\textcolor{black}{Regret-optimal prediction for} the state-space setting}\label{subsec:main_strictly}


In this section, we present the results for the strictly causal scenario of the state-space model in \eqref{eq:def_SS}. Throughout this section, we use the solutions to the Riccati and Lyapunov equation that appear in Section \ref{subsec:riccati}. First, we provide a simple condition for the existence of a predictor of level $\gamma$.
\begin{theorem}[Condition for \textcolor{black}{regret-optimal prediction}]\label{th:strictly_condition}
A $\gamma$-optimal \textcolor{black}{predictor} (the strictly causal scenario) exists if and only if
\begin{align}\label{eq:th_strictly_condition}
    \lambda_{\text{max}}(\overline{Z}_\gamma\Pi)\le 1,
\end{align}
where $\overline{Z}_\gamma$ and $\Pi$ are the solutions to the Lyapunov equations
\begin{align}\label{eq:th_LyaZPI_strictly}
    \Pi&= F_P^\ast \Pi F_P  + H^\ast (I + HPH^\ast)^{-1} H\nn\\
    \overline{Z}_\gamma&= F_P \overline{Z}_\gamma F_P^\ast + (P-U)^\ast L^\ast R_Q^{-1} L (P-U),
\end{align}
and the constants $(F_P,U,R_Q)$ were defined in \eqref{eq:Riccati_P}-\eqref{eq:Lya_U}.
\end{theorem}
The matrices $\overline{Z}_\gamma$ and $\Pi$ have the same dimension as the matrix $F$ in \eqref{eq:def_SS}. Thus, the condition \eqref{eq:th_strictly_condition} boils down to the maximal eigenvalue of a finite-dimensional matrix. The optimal regret can be found using a bisection method on $\gamma$ to guarantee that the maximal eigenvalue in \eqref{eq:th_strictly_condition} is arbitrarily close to $1$. Hereinafter, we assume that the value of $\gamma$ has been optimized so that the $\gamma$-dependent quantities $(W,Q,U,Z_\gamma)$ are fixed.


We proceed to present the solution to the Nehari problem in \eqref{eq:th_reduction_nehari} for the state-space model setting (see remark \ref{remark} on the adaptation to the strictly-causal scenario). Recall that the solution  to the Nehari problem provides the best approximation (in the operator norm sense) for the anticausal part of the operator $\nabla_\gamma K_0\Delta$. In the state-space model this operator can be represented in the frequency domain as $\nabla_\gamma(z) K_0(z)\Delta(z)$, and its anticausal part is denoted by $\overline{T}(z)$ (given below in \eqref{eq:T_strcitly}). The approximation of $\overline{T}(z)$, that is the operator which solves \eqref{eq:th_reduction_nehari} in the strictly-causal scenario is given as follows.
\begin{lemma}[Solution to the Nehari problem - strictly causal]\label{lemma:nehari_SC}
The optimal solution to the Nehari problem (Problem \ref{pro:nehari}) with the anticausal operator $\overline{T}(z)$ (appears in \eqref{eq:T_strcitly} below) is
\begin{align}\label{eq:th_strict_ourNehari}
    \overline{K}_N(z)&= R_Q^{-1/2} L (P-U) \Pi \c{\overline{F}_N} \overline{G}_N,
\end{align}
where
\begin{align}\label{eq:th_strictly_our_nehari_constants}
    \overline{G}_N &=  ( I - F_P \overline{Z}_\gamma F_P^\ast \Pi  )^{-1}F_P \overline{Z}_\gamma H^\ast(I + HPH^\ast)^{-\ast/2}\nn\\
    \overline{F}_N &= F_P - \overline{G}_N (I + HPH^\ast)^{-1/2} H,
    \end{align}
and $\overline{Z},\Pi$ are defined in \eqref{eq:th_LyaZPI_strictly} and are evaluated at the optimal $\gamma$.
\end{lemma}
The solution to the Nehari problem in \eqref{eq:th_strict_ourNehari} corresponds to the operator $\mathcal K'$ that appears in Theorem \ref{th:reduction}. As the solution $\overline{K}_N(z)$ is a strictly causal transfer matrix, the resulted regret-optimal filter in \eqref{eq:th_reduction_filter} will be strictly causal as well. Also, note that the solution in \eqref{eq:th_strict_ourNehari} can be directly transformed to have a state-space form in the time domain. We are now ready to present the regret-optimal predictor.
\begin{theorem}[Regret-optimal predictor]\label{th:strictly_SS}
Given the minimal $\gamma$ that satisfies \eqref{eq:th_strictly_condition}, a regret-optimal \textcolor{black}{predictor} (the strictly-causal scenario) is given by
\begin{align}
    \xi_{i+1} &= \overline{F} \xi_i + \overline{G} y_i\nn\\
    \hat{s}_i &= \overline{H} \xi_i,
\end{align}
where
\begin{align}
    \overline{F}&=
    \begin{pmatrix}
    F_P&0&0\\
    \overline{F}_{2,1}& \overline{F}_N&0\\
    \overline{F}_{3,1}& \overline{F}_{3,2} & F_W
    \end{pmatrix}\nn\\
        \overline{H}&= \begin{pmatrix} L & L(P-U)\Pi& L \end{pmatrix}
    \nn\\
    \overline{G}&= \begin{pmatrix}
    K_P\\
    \overline{G}_N (I + HPH^\ast)^{-1/2}\\
    -(F_Q U + K_Q LP) H^\ast(I + HPH^\ast)^{-1}
    \end{pmatrix},
    \end{align}
    the Riccati constants are defined in \eqref{eq:Riccati_P}-\eqref{eq:Lya_U}, the pair $(\overline{F}_N,\overline{G}_N)$ is defined in \eqref{eq:th_strictly_our_nehari_constants} and the  entries of $\overline{F}$ are
\begin{align*}
    \overline{F}_{2,1} &= - \overline{G}_N (I + HPH^\ast)^{-1/2} H\\
    \overline{F}_{3,1} &=  - (F_Q U + K_Q LP) H^\ast(I + HPH^\ast)^{-1} H \\
    \overline{F}_{3,2} &= K_Q L (P-U) \Pi.
\end{align*}
\end{theorem}
Theorem \ref{th:strictly_SS} concludes the derivation of a simple regret-optimal estimator. In particular, the optimal $\gamma$ can be found using the simple condition in Theorem \ref{th:strictly_condition}, and the estimator can be implemented directly by Theorem~\ref{th:strictly_SS}. Note that the $\gamma-$dependent Riccati equation should be computed throughout the process of determining the minimal $\gamma$, i.e. in Theorem \ref{th:strictly_condition}, but not throughout the filtering process in Theorem \ref{th:strictly_SS}. The internal state of the regret-optimal estimator inherits the finite dimension of the original state space, but it has an increased dimension by a factor of three. Thus, from a computational perspective, the filter requires the same resources as the standard Kalman filter. A simple implementation of the regret-optimal estimator is available in a public Git repository \cite{github_filtering}.

Before concluding this section, it is illuminating to see the relation between the regret-optimal and the Kalman filter via the frequency domain counterpart of Theorem \ref{th:strictly_SS}.
\begin{theorem}[Regret-optimal \textcolor{black}{prediction} in frequency domain]\label{th:strict_freq}
The transfer function of regret-optimal predictor in Theorem \ref{th:strictly_SS} is given by
\begin{align}
    \overline{K}(z)&= \nabla^{-1}_{{\gamma}}(z)[\overline{K}_N(z)+\overline{S}'(z)]\Delta^{-1}(z) + \overline{K}_{H_2}(z),
\end{align}
with
\begin{align}
\overline{S}'(z) &=  - R_Q^{-1/2} L\c{F_Q}\nn\\
 &\ \ \ \cdot (F_Q U-K_Q LP) H^\ast(I + HPH^\ast)^{-\ast/2}\nn\\
 \nabla^{-1}_{{\gamma}}(z) &=  (I+L\c{F_W}K_Q)R_Q^{1/2}\nn\\
\Delta^{-1}(z) &= (I + HPH^\ast)^{-1/2} (I + H (zI-F)^{-1}K_P)^{-1},\nn
\end{align}
$\overline{K}_N(z)$ is defined in \eqref{eq:th_strict_ourNehari} and $\overline{K}_{H_2}(z)$ is the (strictly causal) Kalman filter
\begin{align}
   \overline{K}_{H_2}(z)&= L\c{F_P}K_P.
\end{align}
\end{theorem}
It is interesting that the Kalman filter naturally appears in our derivation of the regret-optimal predictor. This implies that the regret-optimal predictor is a sum of two terms; the first is a Kalman filter which is designed to minimize the Frobenius norm of the operator $T_{\mathcal K}$, while the other term is resulted from a Nehari problem and guarantees that the regret criterion is minimized.

\subsection{\textcolor{black}{Regret-optimal estimation for} the state-space setting}\label{subsec:main_causal}
The structure of this section is similar to the previous section on the strictly causal scenario. We begin with a simple condition for the existence of a regret-optimal estimator.
\begin{theorem}[Condition for regret-optimal estimation]\label{th:condition}
A $\gamma$-optimal causal estimator exists if and only if
\begin{align}\label{eq:th_condition}
    \lambda_{\text{max}}(Z_\gamma\Pi)\le 1,
\end{align}
where $Z_\gamma$ and $\Pi$ are the solutions to the Lyapunov equations
\begin{align}\label{eq:th_condition_LyaZPI}
    \Pi&= F_P^\ast \Pi F_P  + H^\ast (I + HPH^\ast)^{-1} H\nn\\
    Z_\gamma&= F_PZ_\gamma F_P^\ast + F_P(P-U)^\ast L^\ast R_Q^{-1} L (P-U)F_P^\ast.
\end{align}
\end{theorem}


Also here, the solution to the Nehari problem plays an elemental part in the solution. In the causal scenario here, the strictly anticausal part of $\nabla_\gamma(z) K_0(z)\Delta(z)$ is denoted by $T(z)$ (given explicitly in \eqref{eq:T_causal}), and its approximation is given as follows.
\begin{lemma}[Solution to the Nehari problem - causal]\label{lemma:ourNehari}
The optimal solution to the Nehari problem (Problem \ref{pro:nehari}) with the strictly anticausal operator $T(z)$ (Eq. \eqref{eq:T_causal}) is
\begin{align}\label{eq:lemma_nehari_KN}
    K_N(z)&= R_Q^{-1/2} H_N (I + F_N(zI- F_N)^{-1} )G_N,
\end{align}
where
\begin{align}\label{eq:th_our_nehari_constants}
    G_N &=  ( I - F_P Z_\gamma F_P^\ast \Pi  )^{-1}F_P Z_\gamma H^\ast(I + HPH^\ast)^{-\ast/2}\nn\\
    F_N &= F_P - G_N (I + HPH^\ast)^{-1/2} H\nn\\
    H_N&= L (P-U)F_P^\ast \Pi
\end{align}
and the pair $(Z_\gamma,\Pi)$ is defined in \eqref{eq:th_condition_LyaZPI}.
\end{lemma}
The main result for the causal scenario is as follows.

\begin{theorem}[Regret-optimal causal estimation]\label{th:state-space}
Given the optimal threshold $\gamma$, a regret-optimal estimator for the causal scenario is given by
\begin{align}
    \xi_{i+1} &= \tilde{F} \xi_i + \tilde{G} y_i\nn\\
    \hat{s}_i &= \tilde{H} \xi_i + \tilde{J} y_i,
\end{align}
where the matrices are given by
\begin{align}
    \tilde{F}&=
    \begin{pmatrix}
    F_P&0&0\\
    \tilde{F}_{2,1}& F_N&0\\
    \tilde{F}_{3,1}& \tilde{F}_{3,2} & F_W
    \end{pmatrix}\nn\\
    \tilde{H}&= \begin{pmatrix} \tilde{H}_1 &   H_N F_N & L \end{pmatrix}
    \nn\\
    \tilde{G}&= \begin{pmatrix}
    K_P\\ G_N (I + HPH^\ast)^{-1/2}\\ \tilde{G}_3
    \end{pmatrix}\nn\\
    \tilde{J}&= L(P-U)H^\ast(I + HPH^\ast)^{-1} \nn\\
    & \ + H_N G_N (I + HPH^\ast)^{-1/2},
\end{align}
their entries are
\begin{align}
    \tilde{F}_{2,1}&= - G_N (I + HPH^\ast)^{-1/2} H\nn\\
    \tilde{F}_{3,1}&= F_W U H^\ast(I + HPH^\ast)^{-1} H \nn\\
    &\ - K_Q H_N G_N (I + HPH^\ast)^{-1/2} H \nn\\
    \tilde{F}_{3,2}&= K_Q H_N F_N\nn\\
    \tilde{H}_1 &=L- L(P-U)H^\ast(I + HPH^\ast)^{-1} H \nn\\
    &\ - H_N G_N (I + HPH^\ast)^{-1/2} H \nn\\
\tilde{G}_3&= K_Q H_N G_N (I + HPH^\ast)^{-1/2} \nn\\
&\ - F_WU H^\ast(I + HPH^\ast)^{-1},
\end{align}
the Riccati constants are defined in \eqref{eq:Riccati_P}-\eqref{eq:Lya_U} and the variables $(F_N,G_N,H_N)$ are defined in \eqref{eq:th_our_nehari_constants}.
\end{theorem}

\section{Numerical examples}\label{sec:examples}
In this section, numerical experiments are performed in order to evaluate the merits of the regret-optimal filter compared to the traditional $\mathcal H_2$ and $\mathcal H_\infty$ estimators. The experiments are done both in the time domain for particular disturbances, but also in the frequency domain where the behaviour across all energy-bounded disturbances can be illustrated. By \eqref{eq:error_def}, the performance of any linear estimator is governed by the transfer operator $T_{\mathcal K}$ that maps the disturbances $\v{w}$ and $\v{v}$ to the errors sequences $\v{e}$. When the estimator $T_{\mathcal K}$ is time-invariant, it is useful to represent it via its transfer matrix in the $z$-domain
$$
T_K(z) = \left[\begin{array}{cc} L(z) - K(z)H(z) & -K(z) \end{array} \right].
$$
The squared Frobenius norm of $T_K$, which is what the $\mathcal H_2$ estimator minimizes, is given by
\begin{align}\label{eq:def_frob_freq}
\|T_K\|_F^2 = \frac{1}{2\pi}\int_0^{2\pi}\mbox{Tr}\left(T_K^*(e^{j\omega})T_K(e^{j\omega})\right)d\omega,
\end{align}
and the squared operator norm of $T_K$, which is what $\mathcal H_\infty$ estimators minimize, is given by
$$
\|T_K\|^2 = \max_{0\leq\omega\leq 2\pi} \sigma_{\text{max}}\left(T^\ast_K(e^{j\omega})T_K(e^{j\omega})\right),
$$
where $\sigma_{\text{max}}(\cdot)$ denotes the maximal singular value of a matrix.

\begin{figure}[t]
\centering
\psfrag{A}[b][][1]{$\|T_K(e^{j\omega})\|^2$}
\psfrag{B}[t][][1]{Frequency $\omega$}
    \includegraphics[scale=0.4]{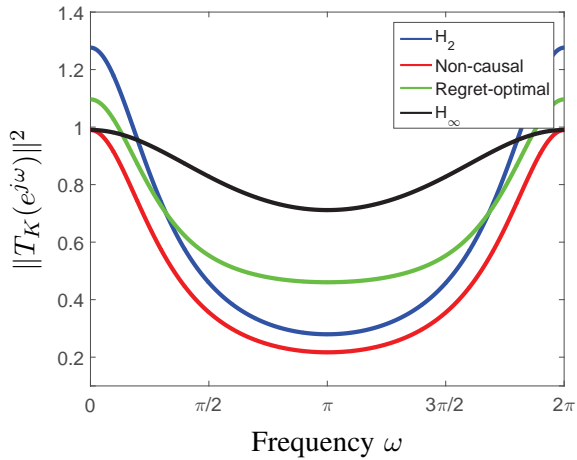}
    \caption{The operator norm as a function of frequency for the scalar system. The curves correspond to the $\mathcal H_2$, $\mathcal H_\infty$ and the regret-optimal estimators, in the causal scenario, and the non-causal estimator. The non-causal estimator achieves the best performance across all frequencies. The $\mathcal H_\infty$ estimator achieves the lowest peak, while the $\mathcal H_2$ estimator attains the smallest area under its curve (i.e., integral) among all causal estimators. The regret-optimal filter attains the best of the two worlds: a significantly lower peak than the $\mathcal H_2$ estimator, and a comparable area under the curve with respect to the optimal $\mathcal H_2$ estimator. Precise comparison of these norms appears in Table \ref{table}.}
    \label{fig:09allestimators}
\end{figure}
\subsection{Scalar system}\label{subsec:scalar}
We start with a simple scalar system to illustrate the results. Consider a stable state-space with $F=0.9$ and $H=L=G=1$. We construct for this state-space the optimal $\mathcal H_2$, $\mathcal H_\infty$, and non-causal estimators, as well as the regret-optimal estimators.
For scalar systems, $T_K(z)$ is a vector so \eqref{eq:def_frob_freq} simplifies as
$$
\|T_K\|_F^2 = \frac{1}{2\pi}\int_0^{2\pi} \left\|T_K(e^{j\omega})\right\|^2 d\omega.
$$
In Fig. \ref{fig:09allestimators}, we plot $\|T_K(e^{j\omega})\|^2$ as a function of frequency for all estimators in the causal scenario. This is insightful as it allows one to assess and compare the performance of the estimators across the full range of input disturbances. Also, the precise norms for the causal and strictly causal scenarios are summarized in Table~\ref{table}.

\begin{table}[b]
\centering
\caption{Estimators performance for the scalar system\label{table}}
\begin{tabular}{c| c| c c c}
                &  & $\|T_K\|_F^2$ & $\|T_K\|^2$ & Regret\\ \hline
  Non-causal & Non-causal estimator & 0.46 & 0.99 &0\rule{0pt}{2.6ex}\\\hline
   & Regret-optimal  & 0.65 & 1.10 & 0.38\rule{0pt}{2.6ex}\\
  Causal &$\mathcal H_2$ estimator & 0.60 & 1.27& 0.70\rule{0pt}{2.6ex}\\ &\rule{0pt}{2.6ex}
 $\mathcal H_\infty$ estimator & 0.94 & 0.99 & 0.71\\ \hline
  & Regret-optimal  & 1.17 & 1.37 & 1.08\rule{0pt}{2.6ex}\\
  Strictly-Causal&$\mathcal H_2$ estimator & 0.99 & 2.21& 2.02\rule{0pt}{2.6ex}\\ &\rule{0pt}{2.6ex}
 $\mathcal H_\infty$ estimator & 1.23 & 1.25 & 1.21
\end{tabular}
\end{table}

\begin{figure}[t!]
\centering
\psfrag{A}[b][][0.8]{$\|T^\ast_K(e^{j\omega})T_K(e^{j\omega})-T^\ast_{K_0}(e^{j\omega})T_{K_0}(e^{j\omega})\|$}
\psfrag{B}[t][][1]{Frequency $\omega$}
    \includegraphics[scale=0.4]{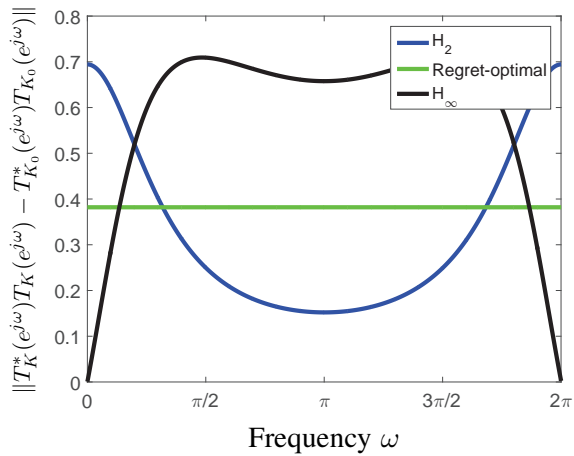}
    \label{fig:09regretall}
    \caption{The regret criterion as a function of frequency. The regret-optimal filter is the one that achieves the lowest peak and it achieves an almost constant regret across all frequencies.}
\end{figure}

It can be observed in Fig. \ref{fig:09allestimators} that the non-causal estimator outperforms the causal estimators across all frequencies. The $\mathcal H_2$ estimator minimizes the Frobenius norm (i.e., the integral above). However, in doing so, it sacrifices the worst-case performance and it has a relatively large peak at low frequencies. The $\mathcal H_\infty$ estimator minimizes the operator norm, i.e., the worst-case performance, which is the peak of its curve (we can see that the $\mathcal H_\infty$ filter attains the minimal $\mathcal H_\infty$ norm since its peak coincides with that of the non-causal estimator). However, in doing so, it sacrifices the average performance measured by the integral. The regret-optimal estimator achieves the best of both worlds by minimizing the largest error deviation from the non-causal filter: its integral is close to that of the $\mathcal H_2$-optimal estimator ($0.6$ vs $0.65$), and it has a peak that significantly improves upon the peak of the $\mathcal H_2$-optimal estimator. We also note in Fig. \ref{fig:09allestimators} that the regret-optimal filter always achieves a smaller cost from either the $\mathcal H_2$ or the $\mathcal H_\infty$ filters. This behaviour will be illustrated in the time-domain simulations of the next example.

It is also illuminating to examine the regret criterion. In Fig. \ref{fig:09regretall}, the regret of the causal estimators is plotted. It can be seen that at low frequencies, the $\mathcal H_\infty$ estimator attains the lowest regret, while at mid-frequencies it is the $\mathcal H_2$ estimator. However, their peak is almost twice than that of the regret-optimal estimator which maintains a balanced, almost constant regret across all frequencies.

\subsection{Tracking example}\label{subsec:tracking}
\begin{figure}[t]
\centering
\psfrag{B}[b][][1]{$\|T_K(e^{j\omega})\|^2$}
\psfrag{A}[t][][1]{Frequency $\omega$}
    \includegraphics[scale=0.4]{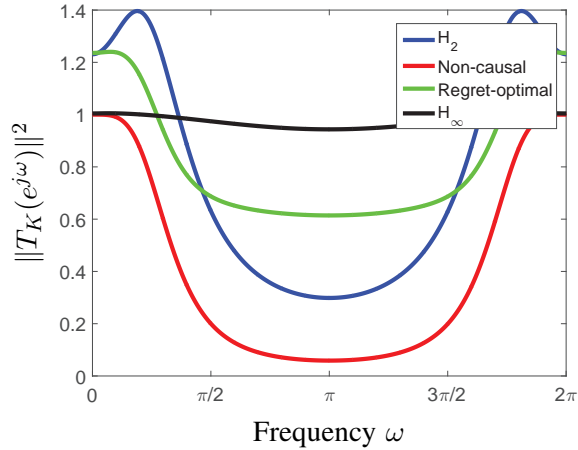}
    \caption{The operator norm of the cost operator as a function of the frequency for the tracking example in Section \ref{subsec:tracking}. Precise comparison for the norms appears in Table \ref{table:tracking}.}
    \label{fig:tracking_freq}
\end{figure}

We study a one-dimensional tracking problem given by
\begin{align}
    \begin{pmatrix}
    x_{i+1} \\ \nu_{i+1}
    \end{pmatrix}&=
    \begin{pmatrix}
    1 & \Delta T\\ 0&1
    \end{pmatrix}    \begin{pmatrix}
    x_{i} \\ \nu_{i}
    \end{pmatrix} +
    \begin{pmatrix}
    0 \\ \Delta T
    \end{pmatrix}a_i\nn\\
    y_i&=      \begin{pmatrix}
    1& 0
    \end{pmatrix}\begin{pmatrix}
    x_{i} \\ \nu_{i}
    \end{pmatrix} + v_i\nn\\
    s_i&=    \begin{pmatrix}
    1& 0
    \end{pmatrix}\begin{pmatrix}
    x_{i} \\ \nu_{i}
    \end{pmatrix},
\end{align}
where $x_i$ corresponds to the position, $\nu_i$ corresponds to velocity and $a_i$ is an exogenous acceleration. The desired signal is the position of the object at the current or the next time step, depending on the availability of the noisy measurements of the noisy position $y_i = x_i + v_i$. In Fig. \ref{fig:tracking_freq}, the frequency response of the various estimators is presented for the causal scenario, and Table \ref{table:tracking} summarizes the precise performance for both the estimation (causal) and the prediction (strictly causal) scenarios.

\begin{figure}[t!]
\centering
\psfrag{B}[b][][1]{Ave. error: $\frac1{t}\sum_{i=1}^t \v{e}^\ast_i \v{e}_i$}
\psfrag{A}[t][][1]{Time $t$}
    \includegraphics[scale=0.2]{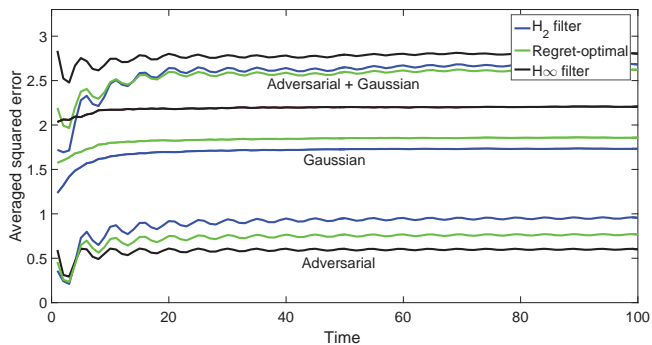}
    \caption{The squared estimation error for the tracking example. There are three triplets of curves where each corresponds to a different disturbances regime: adversarial, Gaussian and their superposition (bottom to top). One can see that the regret-optimal has a close performance to the optimal one in the bottom regimes and outperforms them in their superposition.}
    \label{fig:tracking_time}
\end{figure}

We also simulate the tracking example in time domain to illustrate the performance in different disturbances regimes. In Fig. \ref{fig:tracking_time}, three sets (triplets) of curves describe the filters performance for disturbances with an i.i.d. Gaussian distribution with a variance of $2.25$, adversarial regime (i.e., a worst-case scenario) where the acceleration and the measurement noise are chosen to maximize the transfer function $T_{\mathcal K_{H_2}}$, and the superposition of the two regimes. When the disturbances have a stochastic element, i.e., a Gaussian distribution, the experiment is repeated $10^3$ times and its average is plotted. For disturbances with a Gaussian distribution, the $\mathcal H_2$ filter obtains the smallest estimation error, while for an adversarial noise, the $\mathcal H_\infty$ filter outperforms the others. Interestingly, the regret-optimal filter interpolates between these two filters in both simulations. In fact, for any adversarial disturbances, by Fig.~\ref{fig:tracking_freq}, the regret-optimal filter either interpolates or outperforms the other filters since the green curve lies below the maximum of the blue and black curves. Moreover, there are certain regimes where the superposition of stochastic and adversarial disturbances result in that the regret-optimal filter obtains the smallest error energy, as illustrated in the set of curves at the top of Fig.~\ref{fig:tracking_time}. Also here, this outperforming behavior can be characterized and generalized to other scenarios as long as the Gaussian sequence is added to an adversarial sequence whose frequency lies below $\pi/2$ since the regret-optimal filter has a lower amplitude than the $\mathcal H_2$ filter in this region.

\begin{table}[h]
\centering
\caption{Estimators performance for the tracking problem \label{table:tracking}}
\begin{tabular}{c| c| c c c}
                &  & $\|T_K\|_F^2$ & $\|T_K\|^2$ & Regret\\ \hline
  Non-causal & Non-causal estimator & 0.39 & 1 &0\rule{0pt}{2.6ex}\\\hline
   & Regret-optimal  & 0.82 & 1.24 & 0.65\rule{0pt}{2.6ex}\\
  Causal &$\mathcal H_2$ estimator & 0.77 & 1.40& 1.02\rule{0pt}{2.6ex}\\ &\rule{0pt}{2.6ex}
 $\mathcal H_\infty$ estimator & 0.94 & 0.99 & 0.71\\ \hline
  & Regret-optimal  & 3.82 & 4.00 & 3.80\rule{0pt}{2.6ex}\\
  Strictly-Causal&$\mathcal H_2$ estimator & 3.33 & 6.04 & 5.93\rule{0pt}{2.6ex}\\ &\rule{0pt}{2.6ex}
 $\mathcal H_\infty$ estimator & 3.86 & 3.89 & 3.89
\end{tabular}
\end{table}



\section{Reduction to the Nehari problem (Theorem \ref{th:reduction})}\label{sec:proof_reduction}
This section includes the reduction of the regret problem to the Nehari problem for general linear operators.
\begin{proof}[Proof of Theorem \ref{th:reduction}]
Recall that, for a fixed $\gamma$, the sub-optimal regret problem is given by
\begin{align}\label{eq:proof_reduction_condition_1}
T^\ast_{\mathcal K} T_{\mathcal K} - T_{\mathcal K_0}^\ast T_{\mathcal K_0} \preceq \gamma^2 I.
\end{align}
By the Schur complement, the condition \eqref{eq:proof_reduction_condition_1} is equivalent to
\begin{align}\label{eq:proof_reduction_condition_2}
    &T_{\mathcal K} ( \gamma^{2}I +  T_{\mathcal K_0}^\ast T_{\mathcal K_0})^{-1} T_{\mathcal K}^\ast  \preceq I.
    \end{align}
We use the \emph{Matrix inversion lemma} (in its operator form) to simplify the middle operator on the left hand side of \eqref{eq:proof_reduction_condition_2} as
    \begin{align}\label{eq:proof_operator_MIL}
&\gamma^{-2} I - \gamma^{-2} T_{\mathcal K_0}^\ast (I+ \gamma^{-2} T_{\mathcal K_0}T_{\mathcal K_0}^\ast )^{-1}\gamma^{-2} T_{\mathcal K_0}.
\end{align}

Recall that $\mathcal{K}_0 = \mathcal{L}\mathcal{H}^\ast (I+\mathcal{H}\mathcal{H}^\ast )^{-1}$, so that the cost operator evaluated at the non-causal estimator $\mathcal K_0$ is
\begin{align}
    T_{\mathcal K_0}&= \begin{pmatrix} \mathcal L-\mathcal K_0 \mathcal H & - \mathcal K_0
    \end{pmatrix}\nn\\
    &= \mathcal L(I + \mathcal H^\ast \mathcal H)^{-1} \begin{pmatrix} I & - \mathcal H^\ast
    \end{pmatrix}.
\end{align}
It can now be verified that for any $\mathcal K$,
\begin{align}\label{eq:proof_op_invariantK}
    T_{\mathcal K} T_{\mathcal K_0}^\ast &= \begin{pmatrix} \mathcal L - \mathcal K \mathcal H &- \mathcal K
    \end{pmatrix}\begin{pmatrix} I \\ - \mathcal H
    \end{pmatrix}  (I + \mathcal H^\ast \mathcal H)^{-1} \mathcal L^\ast \nn\\
    &= \mathcal L  (I + \mathcal H^\ast \mathcal H)^{-1} \mathcal L^\ast\nn\\
    &= T_{\mathcal K_0}T_{\mathcal K_0}^\ast.
\end{align}
Using \eqref{eq:proof_operator_MIL} and \eqref{eq:proof_op_invariantK}, we can simplify the condition in \eqref{eq:proof_reduction_condition_2} to
\begin{align}\label{eq:proof_reduction_condition_3}
    T_{\mathcal K} T_{\mathcal K}^\ast & \preceq \gamma^{2}I+ \gamma^{-2} T_{\mathcal K_0}T_{\mathcal K_0}^\ast (I+ \gamma^{-2} T_{\mathcal K_0}T_{\mathcal K_0}^\ast )^{-1} T_{\mathcal K_0}T_{\mathcal K_0}^\ast \nn\\
    &= \gamma^{2}I+ T_{\mathcal K_0}T_{\mathcal K_0}^\ast (\gamma^{2} I+  T_{\mathcal K_0}T_{\mathcal K_0}^\ast )^{-1} T_{\mathcal K_0}T_{\mathcal K_0}^\ast.
\end{align}
The LHS of \eqref{eq:proof_reduction_condition_3} can be written using the completion of the square as
\begin{align}
    T_{\mathcal K}T_{\mathcal K}^\ast&= (\mathcal K - \mathcal K_0)(I + \mathcal H \mathcal H^\ast )(\mathcal K-\mathcal K_0)^\ast + T_{\mathcal K_0}T_{\mathcal K_0}^\ast,\nn
\end{align}
and rearranging the corresponding RHS of \eqref{eq:proof_reduction_condition_3} gives
\begin{align}
    &\gamma^{2}I+ T_{\mathcal K_0}T_{\mathcal K_0}^\ast (\gamma^{2} I+  T_{\mathcal K_0}T_{\mathcal K_0}^\ast )^{-1}T_{\mathcal K_0}T_{\mathcal K_0}^\ast - T_{\mathcal K_0}T_{\mathcal K_0}^\ast\nn\\
    &= \gamma^{2} (I + \gamma^{-2}T_{\mathcal K_0} T_{\mathcal K_0}^\ast )^{-1}.
\end{align}
To conclude, the condition in \eqref{eq:proof_reduction_condition_3} can be written as
\begin{align}\label{eq:proof_reduction_condition_4}
    &(\mathcal K - \mathcal K_0)(I + \mathcal H \mathcal H^\ast )(\mathcal K-\mathcal K_0)^\ast\nn\\
    &\ \ \ \preceq \gamma^{2} (I + \gamma^{-2}\mathcal L(I + \mathcal H^\ast \mathcal H)^{-1} \mathcal L^\ast)^{-1}.
\end{align}
By defining the canonical factorizations
\begin{align}
    \Delta\Delta^\ast &= I + \mathcal H \mathcal H^\ast\nn\\
    \nabla_\gamma ^\ast\nabla_\gamma &= \gamma^{-2}(I + \gamma^{-2}\mathcal L(I + \mathcal H^\ast \mathcal H)^{-1} \mathcal L^\ast),
\end{align}
and applying the Schur complement on \eqref{eq:proof_reduction_condition_4} we get
\begin{align}
    (\mathcal K \Delta- \mathcal K_0\Delta)^\ast \nabla_\gamma ^\ast\nabla_\gamma (\mathcal K\Delta - \mathcal K_0 \Delta)\preceq I
\end{align}
so that
\begin{align}
    (\nabla_\gamma \mathcal K\Delta - \nabla_\gamma \mathcal K_0\Delta)^\ast  (\nabla_\gamma \mathcal K \Delta- \nabla_\gamma \mathcal K_0 \Delta)\preceq I.
\end{align}
Note that $\nabla_\gamma \mathcal K \Delta$ is a causal operator. Now, let $\nabla_\gamma \mathcal K_0 \Delta = \mathcal S + \mathcal T$ where $\mathcal S$ is a causal operator and $\mathcal T$ is a strictly anticausal operator (both operators depend on $\gamma$ implicitly). Then, if
$\mathcal{K}_N$ is a solution to the Nehari problem
\begin{align}
    \|\mathcal{K}_N - \mathcal T\|\le 1,
\end{align}
then a $\gamma$-optimal estimator is given by $\nabla_\gamma^{-1}(\mathcal{K}_N + \mathcal S)\Delta^{-1}$.
\end{proof}

\section{Proof of the state-space setting}\label{sec:proof}
The proofs of the main results for the state-space setting rely on three technical lemmas that are presented in Section \ref{subsec:technical_lemmas}. Based on these technical results, the proofs of the main theorems for the strictly causal and the causal scenarios appear in Section \ref{subsec:proof_str_causal} and Section \ref{subsec:proof_causal}, respectively.

\subsection{Technical lemmas}\label{subsec:technical_lemmas}
By the derivation of Theorem \ref{th:reduction}, there are three technical steps to obtain the Nehari problem:
\begin{enumerate}
    \item Canonical factorization of $\Delta(z)\Delta^\ast(z^{-\ast}) = I+H(z)H^\ast(z^{-\ast}) $.
    \item Canonical factorization of $ \nabla_\gamma(z)\nabla_\gamma^\ast(z^{-\ast}) \\= \gamma^{-2}(I + \gamma^{-2} L(z)(I +  H^\ast(z^{-\ast})  H(z) )^{-1} L^\ast(z^{-\ast}) ).$
    \item Decomposition of $\nabla_\gamma^{-1}(z) L(z)H^\ast(z^{-\ast})\Delta^{-\ast}(z^{-\ast})$.
\end{enumerate}
To preserve the flow of presentation, the proofs of the lemmas appear in Appendix \ref{subapp:proof_lemmas}.

The first step is summarized in the following lemma.
\begin{lemma}[Factorization I]\label{lemma:factor_HHast}
The transfer function $I+H(z)H^\ast(z^{-\ast})$ can be factored as $$\Delta(z)\Delta^\ast(z^{-\ast}) = I+H(z)H^\ast(z^{-\ast}) $$ with
\begin{align}
\Delta(z)&= (I + H (zI-F)^{-1}K_P)  (I + HPH^\ast)^{1/2}
\end{align}
where $(I + HPH^\ast)^{1/2}(I + HPH^\ast)^{\ast/2} = I + HPH^\ast$,  $K_P = FPH^\ast (I + HPH^\ast)^{-1}$ and $P$ is the stabilizing solution to the Riccati equation
\begin{align}
  P &= GG^\ast + FPF^\ast - FPH^\ast (I + HPH^\ast)^{-1} HPF^\ast.\nn
\end{align}
Moreover, the transfer function $\Delta^{-1}(z)$ is bounded on $|z|\ge1$.
\end{lemma}
The first factorization is the one required in the discrete-time Kalman filter and follows from standard arguments that are given for completeness. The second spectral factorization, however, is more involved. In particular, the expression we aim to factor is positive but the order of its causal and anticausal components is in the reversed order in order to perform standard factorization. We overcome this by introducing an additional Riccati equation.
\begin{lemma}[Factorization II]\label{lemma:nablafactor}
For any $\gamma>0$, the factorization
\begin{align}
    &\nabla_\gamma^\ast(z^{-\ast})\nabla_\gamma(z)\nn\\
    &\ = \gamma^{-2}(I + \gamma^{-2}L(z)(I + H^\ast(z^{-\ast}) H(z))^{-1} L^\ast(z^{-\ast}))\nn
\end{align}
holds with
\begin{align}
    \nabla_\gamma(z)&= R_Q^{-1/2}(I-L\c{F_Q}K_Q),
\end{align}
where $R_Q = R_Q^{1/2}R_Q^{\ast/2}$, $Q$ is a solution to the Riccati equation
$$Q = - G R_W^{-1}  G^\ast + F_W Q F_W^\ast  -  K_QR_QK_Q^\ast,$$
and $K_Q= F_W Q L^\ast R_Q^{-1}$ and $R_Q = \gamma^2I + L Q L^\ast$ and the closed-loop system $F_Q = F_W-K_QL$. The constants $(F_W,K_W)$ are obtained from the solution $W$ to the Riccati equation
\begin{align}
    W &= H^\ast H + L_\gamma^\ast L_\gamma + F^\ast W F - K_W^\ast R_W K_W,
\end{align}
and $K_W = R_W^{-1}G^\ast W F$ and $R_W = I + G^\ast W G$ with $R_W = R_W^{\ast/2}R_W^{1/2}$ and $F_W = F - GK_W$.
\end{lemma}
The following lemma is the required decomposition for both scenarios.
\begin{lemma}[Decomposition]\label{lemma:decomposition}
The product of the transfer matrices can be written as
\begin{align}\label{eq:lemma_decompo_product}
 \nabla_\gamma(z)L(z) H^\ast(z^{-\ast})\Delta^{-\ast}(z^{-\ast})   &= \overline{T}(z) + \overline{S}(z)
\end{align}
where $\overline{T}(z)$ is the anticausal transfer matrix
\begin{align}\label{eq:T_strcitly}
        \overline{T}(z) &= R_Q^{-1/2} L (P-U)\nn\\
        &\cdot (I+\a{F_P^\ast}F_P^\ast) H^\ast(I+HPH^\ast)^{-\ast/2},\nn\\
        &= z^{-1} R_Q^{-1/2} L (P-U)\nn\\
        &\cdot \a{F_P^\ast} H^\ast(I+HPH^\ast)^{-\ast/2},
\end{align}
and $\overline{S}(z)$ is the strictly-causal transfer function
\begin{align}\label{eq:S_SC}
     &\overline{S}(z)= \nabla_\gamma(z)L\c{F}F P H^\ast(I \mspace{-3mu} + \mspace{-3mu}HPH^\ast)^{-\ast/2} \\
     & - R_Q^{-1/2} L\c{F_Q}(F_Q U + K_Q LP) H^\ast(I + HPH^\ast)^{-\ast/2}\nn
\end{align}
and $U$ solves $U = K_QL PF_P^\ast + F_QU F_P^\ast$.

Alternatively, the product in \eqref{eq:lemma_decompo_product} can be decomposed for the causal scenario as $T(z) + S(z)$, where $T(z)$ is the strictly anticausal transfer function
\begin{align}\label{eq:T_causal}
    T(z)&= R_Q^{-1/2} L (P\mspace{-3mu}-\mspace{-3mu}U)F_P^\ast\nn\\
    &\ \cdot (z^{-1}I\mspace{-2mu}-\mspace{-2mu}F_P^\ast)^{-1} H^\ast(I \mspace{-3mu}+\mspace{-3mu} HPH^\ast)^{-\ast/2}.
\end{align}
and the causal transfer function is
\begin{align}\label{eq:Scausal}
     &S(z)= \nabla_\gamma(z) L [\c{F}F \mspace{-3mu} + \mspace{-3mu}I]P H^\ast(I \mspace{-3mu} + \mspace{-3mu}HPH^\ast)^{-\ast/2} \nn\\
     & \ - R_Q^{-1/2} L[\c{F_Q}F_Q +I] U H^\ast(I + HPH^\ast)^{-\ast/2}.
\end{align}
\end{lemma}



\subsection{Proof of the strictly causal scenario}\label{subsec:proof_str_causal}
We will jointly prove the results in Section \ref{subsec:main_strictly} concerning the strictly causal scenario.

As mentioned in Remark \ref{remark}, the reduction in the strictly causal scenario directly from Theorem \ref{th:reduction} by following the proof steps and modify \emph{causal} to \emph{strictly causal} and vice versa. The main difference is that in the resulted Nehari problem we have now
\begin{align}
&\min_{\text{s.causal} \ \overline{K}_N(z)}  \| \overline{K}_N(z) - \overline{T}(z)\|.
\end{align}
where $\overline{T}(z)$ is given in \eqref{eq:T_strcitly}. That is, one should approximate an anticausal operator with a \emph{strictly} causal operator. Albeit this formulation is slightly different from the solution in Theorem \ref{th:Nehari_general}, we can transform it into the original Nehari problem as follows
\begin{align}
&\min_{\text{s.causal} \ \overline{K}_N(z)}  \| \overline{K}_N(z) - \overline{T}(z)\|  \nn\\
&\ \ =  \min_{\text{s.causal} \ \overline{K}_N(z)} \| z \overline{K}_N(z) - z \overline{T}(z)\| \nn\\
    &\ \ \stackrel{(a)}=  \min_{ \text{causal} \ K_N(z)} \| K_N(z) - z \overline{T}(z)\|,
\end{align}
where $(a)$ follows from the invertible change of variable $K_N(z) = z\overline{K}_N(z)$. Since the operator $z\overline{T}(z)$ is strictly anticausal, we have a standard Nehari problem. The transfer matrix  of $z\overline{T}(z)$ is also given in the second line of \eqref{eq:T_strcitly}

The solution $\overline{K}_N(z)$ to the Nehari problem in Lemma \ref{lemma:nehari_SC} follows by applying Theorem \ref{th:Nehari_general} with $z\overline{T}(z)$ and multiply the resulted solution with $z^{-1}$ by step $(a)$ above. Furthermore, the condition for the existence of the regret-optimal follows from Theorem \ref{th:Nehari_general} with $z\overline{T}(z)$.

The proof of Theorem \ref{th:strictly_SS} and Theorem \ref{th:strict_freq} follows from the combination of $\overline{K}_N(z)$ in Lemma \ref{lemma:nehari_SC}, the factorizations in Lemma \ref{lemma:factor_HHast} and $\overline{S}(z)$ in Lemma \ref{lemma:nablafactor}. The simplification of the regret-optimal estimator $\nabla^{-1}_\gamma(z) (\overline{S}(z)+\overline{K}_N(z))\Delta^{-1}(z)$ to obtain the explicit time and frequency domains expressions is technical and appears in Appendix \ref{subapp:strictly}.



\subsection{Proof of the causal scenario}\label{subsec:proof_causal}
The proofs of Lemma \ref{lemma:ourNehari} and Theorem \ref{th:condition} on the solution to the Nehari problem and the condition for the estimator existence follow from applying Theorem \ref{th:Nehari_general} with $T(z)$ in \eqref{eq:T_causal}. The proof of Theorem \ref{th:state-space} is technical and appears in Appendix~\ref{subapp:causal}.

\section{Conclusions}\label{sec:conclusion}

We introduced a new criterion for filter design. The regret minimization approach relies on comparing the performance of the designed filter with the performance of a clairvoyant estimator. In the state-space setting, explicit solutions for regret-optimal filtering were obtained for the estimation (the causal scenario) and the prediction (the strictly-causal scenario) regimes. The solutions were obtained by reducing the problem to a Nehari problem. The performance of the regret-optimal filters was compared with the $\mathcal H_2$ and the $\mathcal H_\infty$ estimators for two examples where it was shown to interpolate effectively between the optimal Frobenius and operator norms. 
\appendices

\section{Technical Proofs}
In this section, we prove the technical lemmas in Appendix \ref{subapp:proof_lemmas} and later simplify the estimators for both scenarios.
\subsection{Proofs of the technical lemmas}\label{subapp:proof_lemmas}
The first factorization is a standard spectral canonical  factorization, e.g., \cite[Chapter $12$]{hassibi1999indefinite} \textcolor{black}{but is proved for completeness.}
\begin{proof}[Proof of Lemma \ref{lemma:factor_HHast}]
Consider the term $$I+H(z)H^\ast(z^{-\ast}) = I+H(zI-F)^{-1}GG^\ast(z^{-1}I-F^\ast)^{-1} H^\ast$$ in the matrix form
\begin{align}\label{eq:proof_first_factor}
    & \begin{pmatrix} H (zI-F)^{-1} & I \end{pmatrix}
    \begin{pmatrix} GG^\ast & 0 \\0& I\end{pmatrix}
    \begin{pmatrix} (z^{-1}I-F^\ast)^{-1}H^\ast \\ I \end{pmatrix}
    \nn\\
    &=\textcolor{black}{\begin{pmatrix} H (zI-F)^{-1} & I \end{pmatrix}
    \begin{pmatrix} GG^\ast + FPF^\ast - P& F PH^\ast \\H PF^\ast & I + HPH^\ast\end{pmatrix}}\nn\\
    &\cdot\textcolor{black}{\begin{pmatrix} (z^{-1}I-F^\ast)^{-1}H^\ast \\ I \end{pmatrix}}
    \nn\\
    &= \begin{pmatrix} H (zI-F)^{-1} & I \end{pmatrix}
    \begin{pmatrix} I & \Psi(P) \\ 0 & I \end{pmatrix}
\begin{pmatrix} \Gamma(P) & 0 \\ 0 & I + HPH^\ast\end{pmatrix}\nn\\
&\  \ \cdot \begin{pmatrix} I & 0 \\ \Psi^\ast(P) & I\end{pmatrix}
    \begin{pmatrix} (z^{-1}I-F^\ast)^{-1}H^\ast \\ I \end{pmatrix},
    \end{align}
where the \textcolor{black}{first equality holds for any Hermitian matrix $P$ by direct computation}, and the second equality follows from UDL decomposition with
\begin{align}
  \Psi(P) &\triangleq FPH^\ast (I + HPH^\ast)^{-1}\nn\\
  \Gamma(P) &\triangleq GG^\ast - P + FPF^\ast - FPH^\ast (I + HPH^\ast)^{-1} HPF^\ast.\nn
\end{align}



We now choose $P$ to be the stabilizing solution of the Riccati equation $\Gamma(P)=0$ and denote $K_P\triangleq \Psi(P)$. Then, by \eqref{eq:proof_first_factor}, we can factor the middle matrix as $I + HPH^\ast = (I + HPH^\ast)^{1/2}(I + HPH^\ast)^{\ast/2}$ and conclude
\begin{align}
\Delta(z)&=(I + H (zI-F)^{-1}K_P)  (I + HPH^\ast)^{1/2}.
\end{align}
One can see that the inverse
\begin{align}\label{eq:Dinv}
    \Delta^{-1}(z)&= (I + HPH^\ast)^{-1/2}(I - H (zI-F_P)^{-1}K_P),
\end{align}
is bounded since the closed-loop system $F_P$ is stable.
\end{proof}
The following proof is for the second factorization.
\begin{proof}[Proof of Lemma \ref{lemma:nablafactor}]
Recall that it is required to perform the canonical factorization
\begin{align}\label{eq:proof_secFAC_main}
    &\nabla^\ast(z^{-\ast})\nabla(z)\nn\\
    &\ = \gamma^{-2}(I + \gamma^{-2}L(z)(I + H^\ast(z^{-\ast}) H(z))^{-1} L^\ast(z^{-\ast}))
\end{align}
If we follow a two-step factorization for the inner and outer terms on the RHS of \eqref{eq:proof_secFAC_main}, the transfer matrices will appear in the reversed order. Thus, we should manipulate \eqref{eq:proof_secFAC_main} by inverting it as
\begin{align*}
   &\nabla^{-1}_\gamma(z)\nabla_\gamma^{-\ast}(z^{-\ast})\\
   &= \gamma^2 I \nn\\
   &\ - L(z)(I+H^\ast(z^{-\ast}) H(z) + \gamma^{-2}L^\ast(z^{-\ast}) L(z) )^{-1} L^\ast(z^{-\ast})).
\end{align*}
The factorization will be done in two steps. First, we will factorize the inner expression:
\begin{align}
    \Gamma^\ast(z^{-\ast})\Gamma(z) &= I+H^\ast(z^{-\ast}) H(z) + \gamma^{-2}L^\ast(z^{-\ast}) L(z).
\end{align}
Then, the proof will be completed by factorizing
\begin{align}
    \nabla^{-1}_\gamma(z)\nabla_\gamma^{-\ast}(z^{-\ast})
    &= \gamma^2I - L(z) \Gamma^{-1}(z)\Gamma^{-\ast}(z^{-\ast}) L^\ast(z^{-\ast}).
\end{align}


The inner factorization can be written as
\begin{align}
    &\Gamma^\ast(z^{-\ast})\Gamma(z) \nn\\
    &= I+H^\ast(z^{-\ast}) H(z) + L_\gamma^\ast(z^{-\ast}) L_\gamma(z)\nn\\
    &= I + G^\ast(z^{-1}I-F^\ast)^{-1} (H^\ast H + L_\gamma^\ast L_\gamma)(zI-F)^{-1}G.
\end{align}
Using the same argument as in Lemma \ref{lemma:factor_HHast}, one can show
\begin{align}
    \Gamma(z)&= R_W^{1/2} (I +K_W (zI-F)^{-1}G)
\end{align}
where $W$ is a positive solution for
$$W = H^\ast H + L_\gamma^\ast L_\gamma + F^\ast W F - F^\ast W G (I+G^\ast W G)^{-1}G^\ast W F, $$
and the constants are given by $K_W = (I+G^\ast W G)^{-1}G^\ast W F$ and $R_W=(I + G^\ast W G)$ with $R_W = R_W^{\ast/2}R_W^{1/2}$ and $F_W = F - GK_W$.

Before the second step, we simplify the transfer function
\begin{align}
    &L(z)\Gamma^{-1}(z)\nn\\
    &=L(zI-F)^{-1}[ (I + GK_W (zI-F)^{-1})]^{-1}GR_W^{-1/2}\nn\\
    &= L (zI - F_W)^{-1}GR_W^{-1/2}.
\end{align}
We can now perform the second factorization as
\begin{align}\label{eq:proof_secfactor}
&\nabla^{-1}_\gamma(z)\nabla_\gamma^{-\ast}(z^{-\ast})\\
   &= \gamma^2 I - L(z)\Gamma^{-1}(z)\Gamma^{-\ast}(z^{-\ast}) L^\ast(z^{-\ast})\nn\\
    &= \gamma^2I - L  (zI-F_W)^{-1} G R_W^{-1}  G^\ast (z^{-1}I-F^\ast_W)^{-1}L^\ast\nn\\
    &= (I +L \c{F_W}K_Q)R_Q(I+K_Q^\ast\a{F_W^\ast}L^\ast),\nn
\end{align}
where $Q$ is a solution to the Riccati equation
$$Q = - G R_W^{-1}  G^\ast + F_W Q F_W^\ast  -  F_W QL^\ast (\gamma^2I + L Q L^\ast)^{-1} L Q F_W^\ast, $$
and $K_Q= F_W Q L^\ast R_Q^{-1}$ and $R_Q = \gamma^2I + L Q L^\ast$ and the closed-loop system $F_Q = F_W-K_QL$.

Finally, by \eqref{eq:proof_secfactor}, we can write
\begin{align}
\nabla^{-1}_\gamma(z)&=  (I+L\c{F_W}K_Q)R_Q^{1/2}
\end{align}
with $R_Q = R_Q^{1/2}R_Q^{\ast/2}$ and the proof is completed by taking the inverse
\begin{align}
    \nabla_\gamma(z)&= R_Q^{-1/2}(I+L\c{F_W}K_Q)^{-1}\nn\\
    &= R_Q^{-1/2}(I-L\c{F_Q}K_Q).
\end{align}
\end{proof}
Lastly, the spectral decomposition is proved.
\begin{proof}[Proof of Lemma \ref{lemma:decomposition}]
The decomposition should be performed for the product of the transfer matrices
\begin{align*}
    \nabla_\gamma(z) &= R_Q^{-1/2}(I-L\c{F_Q}K_Q) \\
    L(z)&= L(zI-F)^{-1}G\\
    H^\ast(z^{-\ast})&= G^\ast (z^{-1}I-F^\ast)^{-1}H^\ast\\
    \Delta^{-\ast}(z^{-\ast})&= (I+K_P^\ast(z^{-1}I-F^\ast)^{-1}H^\ast )^{-1}(I + HPH^\ast)^{-\ast/2}.
\end{align*}
First, we simplify $H^\ast(z^{-\ast})\Delta^{-\ast}(z^{-\ast})$ as
\begin{align}
    & G^\ast \a{F_P^\ast} H^\ast(I + HPH^\ast)^{-\ast/2}.
\end{align}

Consider now the product $L(z) H^\ast(z^{-\ast})\Delta^{-\ast}(z^{-\ast})$:
\begin{align}
    &L(z) H^\ast(z^{-\ast})\Delta^{-\ast}(z^{-\ast})\nn\\
    &= L\c{F}G G^\ast \a{F_P^\ast} H^\ast(I + HPH^\ast)^{-\ast/2}\nn\\
    &= L [ \c{F}FP +PF_P^\ast \a{F_P^\ast} +P]\nn\\
    &\cdot H^\ast(I + HPH^\ast)^{-\ast/2},
\end{align}
where in the last equality we use a standard decomposition (e.g., \cite[Lemma $12.3.3$]{hassibi1999indefinite}) detailed as follows. \textcolor{black}{For any matrix $\Sigma$, we can write the product as
\begin{align}\label{eq:proof_first_decomp}
    &\c{F}G G^\ast \a{F_P^\ast}\nn\\
    &= \begin{pmatrix}
    \c{F}&I
    \end{pmatrix}
    \begin{pmatrix}
         GG^\ast - \Sigma + F\Sigma F_P^\ast &F\Sigma\\ \Sigma F_P^\ast&\Sigma
    \end{pmatrix}\nn\\
    &\cdot \begin{pmatrix}
     \a{F_P^\ast}&I
    \end{pmatrix}
\end{align}
By choosing $\Sigma = P$, we note that the first input of the middle matrix satisfies $GG^\ast - P + FPF_P^\ast=0$, and therefore we have the decomposition $\c{F}G G^\ast \a{F_P^\ast} = \c{F}FP +PF_P^\ast \a{F_P^\ast} +P$.}

Now, we can write the causal part of $L(z) H^\ast(z^{-\ast})\Delta^{-\ast}(z^{-\ast})$ as:
\begin{align}
   L [ \c{F}F + I]PH^\ast(I + HPH^\ast)^{-\ast/2}
\end{align}
and its anticausal part as
\begin{align}\label{eq:proof_dec_midanti}
    L  PF_P^\ast \a{F_P^\ast} H^\ast(I + HPH^\ast)^{-\ast/2}.
\end{align}
Note that $\nabla_\gamma(z)$ is causal and, therefore, all
left is to decompose $\nabla_\gamma(z)$ multiplied (from left) with \eqref{eq:proof_dec_midanti}
\begin{align}
    &R_Q^{-1/2}(I-L\c{F_Q}K_Q)L  PF_P^\ast \nn\\
    &\cdot \a{F_P^\ast} H^\ast(I + HPH^\ast)^{-\ast/2}\nn\\
    &= R_Q^{-1/2} L PF_P^\ast \a{F_P^\ast} H^\ast(I + HPH^\ast)^{-\ast/2}\nn\\
    & \ - R_Q^{-1/2} L[\c{F_Q}F_QU + UF_P^\ast \a{F_P^\ast}+U] \nn\\
    &\cdot H^\ast(I + HPH^\ast)^{-\ast/2}\nn\\
    &= R_Q^{-1/2} L (P-U)F_P^\ast \a{F_P^\ast} H^\ast(I + HPH^\ast)^{-\ast/2}\nn\\
    & \ - R_Q^{-1/2} L[\c{F_Q}F_Q +I] U H^\ast(I + HPH^\ast)^{-\ast/2},\nn
\end{align}
where in the first equality, \textcolor{black}{we use the steps in \eqref{eq:proof_first_decomp} to write the decomposition
\begin{align}
    &\c{F_Q}K_Q L P F_P^\ast \a{F_P^\ast}\nn\\&= \c{F_Q}F_QU + UF_P^\ast + U
\end{align}
with $U$ solving the Sylvester equation $U = K_QL PF_P^\ast + F_QU F_P^\ast$.}

To summarize the derivation, the anticausal expression is
\begin{align}\label{eq:Tanticausal}
    T(z)&\mspace{-3mu}= \mspace{-3mu}R_Q^{-1/2} L (P\mspace{-3mu}-\mspace{-3mu}U)F_P^\ast \a{F_P^\ast} \mspace{-3mu} H^\ast\mspace{-3mu}(I \mspace{-3mu}+\mspace{-3mu} HPH^\ast)^{-\ast/2}.
\end{align}
The causal transfer function is
\begin{align}\label{eq:proof_Scausal}
     S(z)&=  \nabla_\gamma(z) L [\c{F}F + I]P H^\ast(I + HPH^\ast)^{-\ast/2}\nn\\
&- R_Q^{-1/2} L[\c{F_Q}F_Q +I] U H^\ast(I + HPH^\ast)^{-\ast/2}
\end{align}
Lastly, it can be verified that $\bar{S}(z)+\bar{T}(z)=S(z) + T(z)$.
\end{proof}

\subsection{Simplification of the strictly-causal estimator}\label{subapp:strictly}
\begin{proof}[Proof of Theorem \ref{th:strictly_SS}]
Recall that by Theorem \ref{th:reduction}, the regret-optimal estimator is given by
\begin{align}
\nabla^{-1}_\gamma(z) (\overline{S}(z)+\overline{K}_N(z))\Delta^{-1}(z).
\end{align}
with the transfer matrices
\begin{align}\label{eq:proof_simplification_sc_all}
        \overline{K}_N(z)&= R_Q^{-1/2} L (P-U) \Pi \c{\overline{F}_N} \overline{G}_N\nn\\
        \Delta^{-1}(z)&= (I + HPH^\ast)^{-1/2} (I + H \c{F}K_P)^{-1}\nn\\
        &= (I + HPH^\ast)^{-1/2} (I - H \c{F_P}K_P)\nn\\
        \nabla^{-1}_\gamma(z)&=  (I+L\c{F_W}K_Q)R_Q^{1/2},
\end{align}
and $\overline{S}(z)$ is given in \eqref{eq:S_SC}.

The first line of $\overline{S}(z)$ in \ref{eq:S_SC} combined with the outer transfer matrices above becomes the standard Kalman filter $L\c{F_P}K_P$. The second line of $\overline{S}(z)$ when multiplied with $\nabla^{-1}_\gamma(z)$ gives
\begin{align}
    & - L\c{F_W}(F_Q U-K_Q LP) H^\ast(I + HPH^\ast)^{-\ast/2},\nn
\end{align}
and multiplying it with $\Delta^{-1}(z)$ from the right gives
\begin{align}\label{eq:prof_sc_simplifiedS}
& - L\c{F_W}(F_Q U + K_Q LP) \nn\\
&\cdot H^\ast(I + HPH^\ast)^{-1}  (I - H \c{F_P}K_P).
\end{align}

The last term that is resulted from the middle term $\overline{K}_N(z)$ cannot be simplified further and thus, we combine this term together with \eqref{eq:prof_sc_simplifiedS} and the Kalman filter to the compact state-space that appears in Theorem \ref{th:strictly_SS}.
\end{proof}

\subsection{Simplification of the causal estimator}\label{subapp:causal}
\begin{proof}[Proof of Theorem \ref{th:state-space}]
Recall that by Theorem \ref{th:reduction}, the regret-optimal estimator is given by
\begin{align}
K(z)&= \nabla^{-1}_\gamma(z) (S(z)+K'(z))\Delta^{-1}(z),
\end{align}
where $\Delta^{-1}(z)$ and $\nabla^{-1}_\gamma(z)$ are given in \eqref{eq:proof_simplification_sc_all}, $K'(z)$ is the solution to the Nehari problem and $S(z)$ is given in \eqref{eq:proof_Scausal}.
First, we take apart the first line of $S(z)$ in \eqref{eq:proof_Scausal} to rewrite
    \begin{align}
K(z)&= \nabla^{-1}_\gamma(z) (S'(z)+K_N(z))\Delta^{-1}(z) + K_{H_2}(z),
\end{align}
where
\begin{align}
     S'(z)&=\mspace{-4mu}  - R_Q^{-1/2}\mspace{-3mu} L[(zI\mspace{-3mu}-\mspace{-3mu}F_Q)^{-1}F_Q\mspace{-3mu} +\mspace{-3mu}I] U H^\ast\mspace{-3mu}(I\mspace{-3mu} +\mspace{-3mu} HPH^\ast)^{-\ast/2},\nn
     \end{align}
and the $\mathcal H_2$ estimator is
\begin{align}\label{eq:proof_simplification_H2_causal}
       K_{H_2}(z)&= LPH^\ast(I + HPH^\ast)^{-1}\nn  \\
       &\ + ( L- LPH^\ast(I + HPH^\ast)^{-1} H)\c{F_P}K_P.
\end{align}
We first simplify the product $\nabla^{-1}_\gamma(z)S'(z)$
\begin{align}
    & - (I+L\c{F_W}K_Q) L U H^\ast(I + HPH^\ast)^{-\ast/2}  \nn\\
    & - L(I-\c{F_Q}K_QL)^{-1} \c{F_Q}F_Q U H^\ast\nn\\
    &\cdot (I + HPH^\ast)^{-\ast/2}  \nn\\
    &= - (I+L\c{F_W}K_Q) L U H^\ast(I + HPH^\ast)^{-\ast/2}  \nn\\
    & - L\c{F_W} F_Q U H^\ast(I + HPH^\ast)^{-\ast/2}\nn\\
    &= - L(I+\c{F_W}F_W) U H^\ast(I + HPH^\ast)^{-\ast/2},\nn
\end{align}
and by multiplying the resulted term with $\Delta^{-1}(z)$ we get
\begin{align}
& \nabla^{-1}_\gamma(z) S'(z) \Delta^{-1}(z)= - L(I+\c{F_W}F_W) U \nn\\
& \cdot H^\ast(I + HPH^\ast)^{-1} (I - H \c{F_P}K_P).
\end{align}

The last term corresponds to the solution of the Nehari problem in Lemma \ref{lemma:ourNehari} and equals
\begin{align}\label{eq:proof_simplifcation_third}
    &\nabla^{-1}(z)R_Q^{-1/2} H_N (I + F_N(zI- F_N)^{-1} )G_N \Delta^{-1}(z).
\end{align}
It can now be directly verified that the sum of \eqref{eq:proof_simplification_H2_causal}-\eqref{eq:proof_simplifcation_third} can be compactly expressed as $K(z)= \tilde{H}\c{\tilde{F}}\tilde{G} + \tilde{J}$ with the matrices given in Theorem \ref{th:state-space}.

\end{proof}

\bibliography{ref}
\bibliographystyle{IEEEtran}
\end{document}